\documentclass[12 pt]{article}
\usepackage{amsmath,amsfonts,fancybox,multicol,lscape}
\usepackage{color}
\usepackage{graphicx}
\usepackage{tikz}
\usetikzlibrary{arrows,automata,positioning}
\usetikzlibrary{quotes,angles}

\newtheorem{lemma}{Lemma}
\newtheorem{corollary}{Corollary}

\newtheorem{theorem}{Theorem}

\begin{document}

\title{A Degree Six Inequality on 

Convex Quadrilaterals}

\author{Mazen Bou  Khuzam\\
Dept. of Mathematics and Natural Science \\
The American University of Iraq at Sulaimani  \\
Email: mazen.boukhuzam@auis.edu.krd}
\date{}
\maketitle

\begin{abstract}
We prove a degree-six inequality on convex quadrilaterials. This inequality originated from work on the Atiyah-Sutcliffe conjectures on configurations of points in $\Bbb{R}^{3}$.
\end{abstract}

\noindent $\bf{Key words}$: Atiyah determinant ; Atiyah-Sutcliffe conjectures

\section{Introduction}

\let\thefootnote\relax\footnote{ I would like to thank Ms. Shams Watha for her help in writing numerous codes to verify some of the inequalities associated with this work}

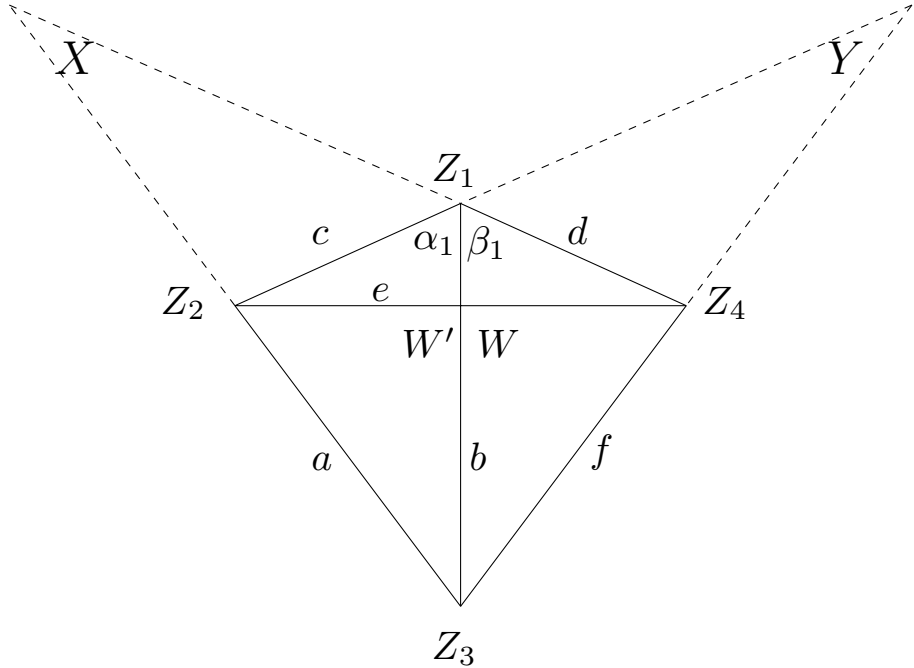
\begin{figure}[htb!]

\centering

\begin{tikzpicture}
[scale=2, every node/.style={anchor=center, scale=1.5}]
  \draw
    (-2,0) coordinate (a) node[right] {}
    (4,0) coordinate (b) node[left] {}
    (1,-1.32) coordinate (c) node[above right] {}
    (0.7,-1.32) coordinate (z4) node[above right] {{\footnotesize $Z_1$}}
    (-0.5,-2) coordinate (d) node[above right] {}
    (-1.1,-2.2) coordinate (z1) node[above right] {{\footnotesize $Z_2$}}
    (2.5,-2) coordinate (e) node[above right] {}
     (0.5,-2.45) coordinate (z1) node[above right] {{\footnotesize $W'$}}
    (2.5,-2) coordinate (e) node[above right] {}
    (1.0,-2.45) coordinate (z1) node[above right] {{\footnotesize $W$}}
    (2.5,-2) coordinate (e) node[above right] {}
  (-0.1,-1.7) coordinate (z1) node[above right] {{\footnotesize $c$}}
    (2.5,-2) coordinate (e) node[above right] {}
  (1.6,-1.7) coordinate (z1) node[above right] {{\footnotesize $d$}}
    (2.5,-2) coordinate (e) node[above right] {}
  (1.75,-3.2) coordinate (z1) node[above right] {{\footnotesize $f$}}
    (2.5,-2) coordinate (e) node[above right] {}
(0.95,-3.2) coordinate (z1) node[above right] {{\footnotesize $b$}}
    (2.5,-2) coordinate (e) node[above right] {}
  (-0.1,-3.2) coordinate (z1) node[above right] {{\footnotesize $a$}}
    (2.5,-2) coordinate (e) node[above right] {}
 (0.3,-2.08) coordinate (z1) node[above right] {{\footnotesize $e$}}
    (2.5,-2) coordinate (e) node[above right] {}
    (2.5,-2.2) coordinate (z3) node[above right] {{\footnotesize $Z_4$}}
    (1,-4) coordinate (f) node[above right] {}
    (0.7,-4.5) coordinate (z2) node[above right] {{\footnotesize $Z_3$}};  
 \draw[dashed]  (-2,0) -- (1,-1.32);
 \draw   (1,-1.32) -- (2.5,-2);
   \draw[dashed] (-0.5,-2) -- (-2,0);
\draw (-0.5, -2) -- (1,-4);
\draw[dashed] (4,0) -- (2.5,-2);  
\draw (2.5,-2) -- (1,-4);
\draw[dashed]    (4,0) -- (1,-1.32);
\draw   (1,-1.32) -- (-0.5,-2)
    (-0.5,-2) -- (2.5,-2)
    (1,-1.32) -- ((1,-4)
    pic["$X$", angle eccentricity=1.5, angle radius=0.5cm]
    {angle=f--a--e}
    pic["$Y$", angle eccentricity=1.5, angle radius=0.5cm]
    {angle=d--b--f}
    pic["{\footnotesize $\alpha_1$}", angle eccentricity=1.7, angle radius=0.25cm]
    {angle=d--c--f}
    pic["{\footnotesize $\beta_1$}", angle eccentricity=1.7, angle radius=0.25cm]
    {angle=f--c--e};
\end{tikzpicture}
 
\caption{Convex Quadrilateral}
\label{fig1}
\end{figure}

The inequality we introduce  in this paper has its origin from the Atiyah determinant. M. F. Atiyah constructed this determinant  in \cite{A1} (see also related work in \cite{A2}) in an attempt to answer a geometric question asked by  M. V. Berry and J. M. Robbins in \cite{BR} based on their work on the Spin Statistics Theorem. He conjectured that the determinant should be non-zero in order to  finalize one of his original constructions. Together with P. Sutcliffe, he added two new conjectures in \cite{AS}  implying the original conjecture. However, proving these conjectures turned out to be notoriously difficult, as described by him. The linear case and the case of three points were rather straightforward, but the problem started to become increasingly resistent when we reach the case of four points. M. Eastwood and P. Norbury proved the first conjecture in \cite{EN} for the case of four points, and other proofs for special cases have appeared in \cite{D} , \cite{MP}, and \cite{M}. Nevertheless, the main reason or geometric reason behind the difficulty of the $n=4$ case and the general case remained in question. In this paper, we introduce a new geometric inequality obtained from work on the Atiyah determinant. This degree-six inequality seems to lie at the heart of the problem and  will be used in a subsequent paper to prove the second  Atiyah-Sutcliffe conjecture on convex quadrilaterals. We hope that this work would help in finding similar  inequalities for $n \geq 5$ in order to solve the general conjecture.


To present this  inequality, we start by assuming we are working in the complex plane and our four points are located at  the pairwise distinct points $z_{1}$, $z_{2}$, $z_{3}$, and $z_{4}$. The case of  a convex quadrilateral is shown in figure 1. In the triangle $z_{1}z_{2}z_{3}$, the lengths of the sides are assigned as follows: The side facing the vertex at $z_{1}$ has length $a$,  the side facing the vertex at $z_{2}$ has length $b$, and the side facing the vertex at $z_{3}$ has length $c$. The other three sides have lengths   $|z_{4}z_{1}|=d$, $|z_{4}z_{2}|=e$, and $|z_{4}z_{3}|=f$. We define the following degree-six expression corresponding to the side $z_{1}z_{2}$ by

\[ E_{12} = f \, A_{123} \, A_{124} \, (a+b+e+d-2c) \]

\noindent  This expression includes the area of the triangle $z_{1}z_{2}z_{3}$, denoted by $A_{123}$, multiplied by the area of the triangle $z_{1}z_{2}z_{4}$. These are the two triangles sharing the side $z_{1}z_{2}$. This is multiplied with $(a+b-c)+(e+d-c)$ which is the sum of the two triangle inequalities at the side  $z_{1}z_{2}$ which has length $c$. Finally, the whole product is multiplied by $f$, the length of the free (unused) side. It is obvious that $E_{12}$ is of degree six and is positive by construction. One can define five similar quantities, one for each of the remaining five sides of the quadrilateral. For example: $ E_{13}= e \,  A_{123} \, A_{134} \, (c+a+d+f-2b)$. Using these six expressions we introduce and prove the following inequality

\begin{theorem}

 For a convex quadrilateral, the following inequality is satisfied:

\[ E_{12}+ E_{23} + E_{34}+E_{41} \geq  E_{13} +  E_{24}   \tag{1} \]

\end{theorem}

In section 2 we simplify the expression of the multiplicity-one terms. In section 3, we simplify the multiplicity-two terms and then we use these two simplifications to prove the inequality in section 4.

\section{The  Multiplicity-One Terms}

Each of the $E_{ij}$'s consists of 5 terms which means that our inequality comprise of 30 terms. We start by displaying these thirty terms in order to be able to select from them based on symmetric considerations that will be used in the proof of the inequality.  Here is the full inequality: 

\vspace{3 mm}

$\displaystyle{    f  \, A_{123} \, A_{124} \, (d+e+a+b-2c) +d \, A_{123} \, A_{234} \, (c+b+e+f-2a) }$

$\displaystyle{   +  c \, A_{134} \, A_{234} \, (d+b+e+a-2f)  +  a \, A_{124} \, A_{134} \, (c+e+b+f-2d) } $

$\displaystyle{ -  e \,  A_{123} \, A_{134} \, (c+a+d+f-2b)      -  b \, A_{124} \, A_{234} \, (c+d+a+f-2e)      \geq 0 }$ \begin{equation*}
 \tag{2}
\end{equation*}

The first important note about inequality (2) is that the terms of multipicity one happen to have a natural decomposition into three disjoint classes: The $X$-terms, the $Y$-terms, and the $W$-terms, where $X$, $Y$, and $W$ are the angles shown in figure 1. We will start with the $X$-terms first: These are multiplicity-one terms in (2) that do not have an $a$ or $d$ factor. We call these the $X$-terms since $a$ and $d$ are the side-lengths of the angle $X$. Accordingly, the $X$-terms are:

\vspace{3 mm}

\noindent $\displaystyle{f(e+b) A_{123}  A_{124} + c(b+e)  A_{134}  A_{234}  - e(c+f)  A_{123}  A_{134}  - b (c+f)  A_{124}  A_{234} }$

\vspace{2 mm}

\noindent Rearranging these terms gives us:

\vspace{3 mm}

\noindent $\displaystyle{ bc \, A_{234} ( A_{134} - A_{124} ) + ce \, A_{134} \, ( A_{234} - A_{123} ) - bf A_{124} \, (   A_{234} - A_{123} )}$

\noindent $\displaystyle{ -ef \, A_{123} \, (A_{134} - A_{124} ) }$

\vspace{2 mm}

\noindent We note that the last factor in each of the above four terms is actually the same. It is easy to see that this factor  is equal to  $\displaystyle{\frac{1}{2} ad \sin X }$. Substituting this back in the last expression, we see that the $X$-terms can be written as: 

\[ \left(  bc \, A_{234} + ce \, A_{134} -  bf A_{124}  - ef \, A_{123} \right)  (\frac{1}{2} a d \sin X ) \]

\noindent Replacing the areas of the triangles using the formulae $A_{234}= \frac{1}{2} ef \sin \beta _{4} $, $A_{134}= \frac{1}{2} bf \sin \alpha _{3} $, $A_{124}= \frac{1}{2} ce \sin \beta _{2} $, and $A_{123}= \frac{1}{2} bc \sin \alpha _{1} $ allows us to decompose the $X$-terms as a product of a scalar quantity with an angular one:

\[ \frac{1}{4} (abcdef) ( \sin \beta _{4} + \sin \alpha _{3} - \sin \beta _{2} - \sin \alpha _{1}) \sin X \]

 \noindent The angular part can be simplified further by using that $\displaystyle{\sin \alpha + \sin \beta = }$

\noindent $2 \sin ( \frac{\alpha + \beta}{2}) \, \cos ( \frac{\alpha - \beta}{2}) $  on the first two terms and the last two terms to obtain

\[ \frac{1}{2} (abcdef)  \left(   \sin ( \frac{\alpha _{3}  +  \beta _{4}}{2}) \cos ( \frac{\alpha _{3}  -  \beta _{4}}{2}) -  \sin ( \frac{ \alpha _{1} + \beta _{2}}{2}) \cos ( \frac{ \alpha _{1} - \beta _{2}}{2}) \right) \sin X \]

\noindent Since $W' = \alpha _{3}  +   \beta _{4} =   \alpha _{1} +  \beta _{2}$, the answer becomes:

\[ \frac{1}{2} (abcdef)   \left( \cos ( \frac{\alpha _{3}  -  \beta _{4}}{2}) - \cos ( \frac{ \alpha _{1} - \beta _{2}}{2}) \right)  \sin \frac{W'}{2} \sin X \]

\noindent Finally, using $\displaystyle{\cos \alpha - \cos \beta = -2 \sin ( \frac{\alpha + \beta}{2}) \, \sin ( \frac{\alpha - \beta}{2}) }$, and the fact that

\noindent $\alpha _{3}  - \beta _{2} = \alpha _{1} - \beta _{4}$ and $ \alpha _{3} - \alpha _{1} =  \beta _{2}- \beta _{4} = -Y$, the $X$-terms take the factorized form:

\[  (abcdef) \sin X  \sin \frac{W'}{2}  \sin \frac{Y}{2} \sin ( \frac{\alpha _{1}-\beta _{4}}{2}) \]

\noindent Similar factorrization for the $Y$-terms and $W$-terms can be obtained in exactly the same way and are summarized in the following lemma:

\vspace{3 mm}

\begin{lemma}: The multiplicity-one terms of our inequality can be simplified as follows:

The $X$-terms add up to: $\displaystyle{  (abcdef) \sin X  \sin \frac{W'}{2}  \sin \frac{Y}{2} \sin ( \frac{\alpha _{1}-\beta _{4}}{2}) }$

\vspace{3 mm}

The $Y$-terms add up to: $\displaystyle{  (abcdef) \sin Y  \sin \frac{W}{2}  \sin \frac{X}{2} \sin ( \frac{\beta _{1}-\alpha _{2}}{2} )}$

\vspace{3 mm}

The $W$-terms add up to: $\displaystyle{  (abcdef) \sin W  \cos \frac{X}{2}  \cos \frac{Y}{2} \sin ( \frac{\gamma _{1}+\gamma _{3}}{2} )}$

\end{lemma}

\section{The Multiplicity-Two Terms}

We will see that these multiplicity-two terms have a natural dichotomy requiring different ways for simplifying their expressions. First, let us collect them from the expression in (2) and display them in front of us:

\vspace{2 mm}

\noindent $\displaystyle{ -2ad A_{123}A_{234} - 2ad A_{124}A_{134} -2cf A_{124} A_{123} -2cf A_{134}A_{234} }$

\noindent $\displaystyle{+ 2be A_{123}A_{134} + 2be A_{124}A_{234}}$

\vspace{2 mm}

\noindent Similar to what we did for multiplicity-one terms, we write this expression as a scalar part multiplied with an angular part as follows:

\vspace{2 mm}

\noindent $\displaystyle{ \frac{1}{2} (abcdef) ( -\sin \alpha _{1} \sin \beta _{4}  -\sin \alpha _{3} \sin \beta _{2} -\sin \alpha _{4} \sin \beta _{3} -\sin \alpha _{2} \sin \beta _{1}  }$

\noindent \hspace{2 cm} $\displaystyle{+\sin \gamma _{1} \sin \gamma _{3} -\sin \gamma _{2} \sin \gamma _{4} ) }$

\vspace{2 mm}

\noindent We linearize each term using $\displaystyle{\sin \alpha \sin \beta = \frac{1}{2} (\cos (\alpha - \beta ) - \cos (\alpha + \beta ))}$ and regroup their angular parts as follows:

\vspace{2 mm}

\noindent The first part of the angular expression is

\vspace{2 mm}

 \noindent $\displaystyle{P_{1} = \frac{1}{4}  (\cos (\alpha _{1} + \beta _{4}) +\cos (\alpha _{3} + \beta _{2}) +\cos (\alpha _{4} +\beta _{3}) + \cos (\alpha _{2}+ \beta _{1})}$

$\displaystyle{ +\cos (\gamma _{1}-\gamma _{3}) + \cos (\gamma _{2} - \gamma _{4}) )}$
\vspace{2 mm}

\noindent The second part of the angular  expression is

\vspace{2 mm}

 \noindent $\displaystyle{P_{2} = \frac{1}{4}  (-\cos (\alpha _{1} - \beta _{4}) -\cos (\alpha _{3} - \beta _{2}) -\cos (\alpha _{4} -\beta _{3}) - \cos (\alpha _{2}- \beta _{1})}$

$\displaystyle{ -\cos (\gamma _{1}+\gamma _{3}) - \cos (\gamma _{2} + \gamma _{4}) )}$

\vspace{2 mm}

\noindent  We notice that all angles appearing in the expression of $P_{1}$ can be expressed in terms of $X$, $Y$, and $W$ (or $W' = \pi -W$). Namely, $\alpha _{1} + \beta _{4} = W' + Y$, $\alpha _{3} + \beta _{2} = W' - Y$, $\alpha _{2} + \beta _{1} = W + X$, $\alpha _{4} + \beta _{3} = W - X$, $\gamma _{1} - \gamma _{3} =  X + Y$, and $\gamma _{2} - \gamma _{4} = X - Y$. Accordingly, $P_{1}$ can be simplified as follows:

\vspace{2 mm}

\noindent $\displaystyle{P_{1}= \frac{1}{2} (  \cos X \cos W +  \cos Y \cos W' + \cos X \cos Y ) }$

\vspace{2 mm}

\noindent This can be simplified further by using the half-angle formula to produce

\vspace{2 mm}

\noindent $\displaystyle{P_{1}= \frac{1}{2}  (1-2\sin ^{2} \frac{X}{2})(1-2\sin ^{2} \frac{W}{2}) +  \frac{1}{2}(1-2\sin ^{2} \frac{Y}{2}) (1-2\sin ^{2} \frac{W'}{2}) }$

\hspace{3 mm} $ +  \frac{1}{2} (1-2\sin ^{2} \frac{X}{2}) (1-2\sin ^{2} \frac{Y}{2}) $

\vspace{2 mm}

\noindent Using $\displaystyle{\sin ^{2} \frac{W}{2} + \sin ^{2} \frac{W'}{2} =1}$ we can simplify $P_{1}$ to become

\vspace{2 mm}

\noindent $\displaystyle{P_{1}= \frac{1}{2} -2 \sin ^{2} \frac{X}{2} \cos ^{2} \frac{W}{2} -2\sin ^{2} \frac{Y}{2} \cos ^{2} \frac{W'}{2} + 2 \sin ^{2} \frac{X}{2}   \sin ^{2} \frac{Y}{2}  }$

\vspace{2 mm}

\noindent Inserting $\displaystyle{  \sin ^{2} \frac{Y}{2}  + \cos ^{2} \frac{Y}{2}  }$ as a factor in the second term, $\displaystyle{  \sin ^{2} \frac{X}{2} +  \cos ^{2} \frac{X}{2}  }$ as a factor in the third term, and $\displaystyle{  \sin ^{2} \frac{W}{2}  + \cos ^{2} \frac{W}{2}  }$ as a factor in the fourth term allows us to write $P_{1}$ as

\vspace{1 cm}

\[ P_{1} =  \frac{1}{2} - 2  \sin ^{2} \frac{X}{2} \cos ^{2} \frac{W}{2} \cos ^{2} \frac{Y}{2}   - 2  \cos ^{2} \frac{X}{2} \cos ^{2} \frac{W'}{2}  \sin ^{2} \frac{Y}{2} \]

\noindent As for $P_{2}$, we notice that $\alpha _{1} - \beta _{4} = \alpha _{3} - \beta _{2}$, $\alpha _{4} -\beta _{3} =\alpha _{2}- \beta _{1}$, and $\gamma _{1}+\gamma _{3} = 2\pi  - (\gamma _{2} + \gamma _{4})$ which  allows us to write it as 

\vspace{2 mm}

\noindent $\displaystyle{ P_{2} = \frac{1}{2} \left( - \cos (\beta _{4} - \alpha _{1}) - \cos ( \alpha _{2} - \beta _{1} ) - \cos ( \gamma _{3} + \gamma _{1}) \right) }$. By realizing that the sum of the angles in this last expression is $\pi$ and using the identity 

\[ \cos \alpha + \cos \beta + \cos \gamma = 1 + 4 \sin (\frac{\alpha}{2}) \sin (\frac{\beta}{2}) \sin (\frac{\gamma}{2}) \hspace{6 mm} ; \hspace{3 mm} \alpha + \beta + \gamma = \pi \]

\begin{equation*}
 \tag{3}
\end{equation*}

\vspace{2 mm}

\noindent we  can see that $\displaystyle{P_{2}= -\frac{1}{2} - 2 \sin (\frac{\alpha _{1} - \beta _{4}}{2}) \sin (\frac{\beta _{1} - \alpha _{2} }{2}) \sin (\frac{\gamma _{1} + \gamma _{3} }{2})}$

\vspace{2 mm}

\noindent We add $P_{1}$ and $P_{2}$ and summarize these computations in the following lemma

\vspace{2 mm}

\begin{lemma}: The multiplicity-two terms of our inequality can be writen as $ (abcdef)$ times

\[ - 2  \sin ^{2} \frac{X}{2} \cos ^{2} \frac{W}{2} \cos ^{2} \frac{Y}{2}   - 2  \cos ^{2} \frac{X}{2} \cos ^{2} \frac{W'}{2}  \sin ^{2} \frac{Y}{2}  \]

\[ - 2 \sin (\frac{ \alpha _{1} - \beta _{4}}{2}) \sin (\frac{  \beta _{1} - \alpha _{2}}{2}) \sin (\frac{ \gamma _{1} + \gamma _{3}}{2}) \]

\end{lemma}

\section{Proof of The Inequality}

Now that we have simplified each group of terms in the previous two sections, we turn to the proof of inequality (1). We first prove the that the sum of the $X$, $Y$, $W$-terms, and $P_{1}-\frac{1}{2}$ (which is the first line of Lemma 2) is non-negative. In particular, using the simplifications obtained in the previous two sections, we prove the following theorem:

\vspace{2 mm}

\begin{theorem}: For the convex quadrilateral presented in figure 1 where $\gamma _{2} + \gamma _{3} \leq \pi$ 
 and $\gamma _{3} + \gamma _{4} \leq \pi$, we have 

\vspace{4 mm}

$\begin{array}{cccc}
\sin X  \sin \frac{W'}{2}  \sin \frac{Y}{2} \sin ( \frac{\alpha _{1}-\beta _{4}}{2} )      &                &       & \\
+   \sin Y  \sin \frac{W}{2}  \sin \frac{X}{2} \sin ( \frac{\beta _{1}-\alpha _{2}}{2} )   &  \begin{array}{c}
- 2  \sin ^{2} \frac{X}{2} \cos ^{2} \frac{W}{2} \cos ^{2} \frac{Y}{2} \\
 \\
 - 2  \cos ^{2} \frac{X}{2} \cos ^{2} \frac{W'}{2}  \sin ^{2} \frac{Y}{2}  
\end{array} & \geq & 0 \\
+   \sin W  \cos \frac{X}{2}  \cos \frac{Y}{2} \sin ( \frac{\gamma _{1}+\gamma _{3}}{2} ) &        &     & 
\end{array}$ 

 \begin{equation*} \tag{4} \end{equation*}

\end{theorem}

\vspace{2 mm}

\noindent The first column in this inequality is nothing but the simplified expressions of the sum of the $X$, $Y$, and $W$-terms, whereas the second column is the first line in Lemma 2 (which is $P_{1}-\frac{1}{2}$). To prove Theorem 2, we first need the following corollary:

\begin{corollary}:  On the convex quadrilateral depicted in figure 1, we have the following 3 inequalities:

\[ \left| \sin \left( \frac{ \alpha _{1} -\beta _{4}}{2} \right) \right| \leq \sin \left( \frac{  W' -Y}{2} \right) \] 

\[ \left| \sin \left( \frac{ \beta _{1} -\alpha _{2}}{2} \right) \right| \leq \sin \left( \frac{  W - X}{2} \right)\]

\[  \sin \left( \frac{ \gamma _{1} +\gamma _{3}}{2} \right)  \geq   \sin \left( \frac{ X +Y}{2}  \right)   \] 
\end{corollary}

\vspace{4 mm}

\noindent Proof: By using $ \sin ( \frac{ \alpha _{3} -\beta _{2}}{2} ) +  \sin ( \frac{ \alpha _{3} +\beta _{2}}{2} ) = 2 \sin  \frac{ \alpha _{3}}{2} \cos \frac{\beta _{2}}{2} $ and noticing that

\vspace{2 mm}

\noindent $ \alpha _{3} -\beta _{2} = \alpha _{1} -\beta _{4}$ and $\alpha _{3} +\beta _{2} = W' - Y$, we conclude:

\vspace{3 mm}

$ \sin ( \frac{ \alpha _{1} -\beta _{4}}{2} ) =  \sin ( \frac{ Y - W'}{2} ) + 2 \sin  \frac{ \alpha _{3}}{2} \cos \frac{\beta _{2}}{2} \geq \sin ( \frac{ Y - W'}{2} ) $

\vspace{2 mm}

\noindent In addition, $\sin (\frac{\alpha _{3} - \beta _{2}}{2} ) \leq \sin (\frac{\alpha _{3} + \beta _{2}}{2} ) $ implies that $\sin (\frac{\alpha _{1} - \beta _{4}}{2} ) \leq \sin (\frac{W' - Y}{2} ) $ which concludes the first inequality in this corollary.

\vspace{2 mm}

\noindent The other two  inequalities can be proved in a similar way by using 

\vspace{3 mm}

$ \sin ( \frac{ \beta _{3} -\alpha _{4}}{2} ) +  \sin ( \frac{ \beta _{3} +\alpha _{4}}{2} ) =  2 \sin  \frac{ \beta _{3}}{2} \cos \frac{\alpha _{4}}{2} $ for the second one and 

\vspace{3 mm}

$ \sin ( \frac{ \gamma _{1} +\gamma _{3}}{2} ) -  \sin ( \frac{ \gamma _{1} -\gamma _{3}}{2} ) = 2 \cos  \frac{ \gamma _{1}}{2} \sin \frac{\gamma _{3}}{2} $ for the third one.

\vspace{3 mm}

\noindent  We also have to use that  $ \beta _{3} -\alpha _{4} = \beta _{1} -\alpha _{2}$ and $\beta _{3} +\alpha _{4} = W - X$ for the second inequality and $ \gamma _{1} -\gamma _{3} = X+Y$ for the last one.

\vspace{5 mm}

\noindent $\bf{Proof \, \,  of \, \, Theorem \, 2}$: Based on Corollary 1, the left hand side of (4)  is not less than:

\vspace{2 mm}

$\begin{array}{ccc}  
\sin X  \sin \frac{W'}{2}  \sin \frac{Y}{2} \sin ( \frac{Y-W'}{2} )  &  & \\

+  \sin Y  \sin \frac{W}{2}  \sin \frac{X}{2} \sin ( \frac{X-W}{2} ) & \begin{array}{c}
 - 2  \sin ^{2} \frac{X}{2} \cos ^{2} \frac{W}{2} \cos ^{2} \frac{Y}{2} \\
  - 2  \cos ^{2} \frac{X}{2} \cos ^{2} \frac{W'}{2}  \sin ^{2} \frac{Y}{2} 
\end{array}
& \\
+  \sin W  \cos \frac{X}{2}  \cos \frac{Y}{2} \sin ( \frac{X+Y}{2} ) & & 
\end{array}$

\vspace{2 mm}

\noindent Expanding  $\displaystyle{\sin  (X+Y)/2  =\sin X/2 \cos Y/2 + \cos X/2 \sin Y/2}$ in the third line of the first column and combining the terms with the ones in the second column yields:

\vspace{2 mm}

$\begin{array}{cc}
 2  \sin \frac{X}{2} \cos \frac{X}{2}  \sin \frac{W'}{2}  \sin \frac{Y}{2} \sin ( \frac{Y-W'}{2} ) &   \\
 + 2 \sin \frac{Y}{2}  \cos \frac{Y}{2} \sin \frac{W}{2}  \sin \frac{X}{2}  \sin ( \frac{X-W}{2} ) & \\
 + 2 \cos \frac{W}{2} \sin \frac{X}{2} \cos ^{2} \frac{Y}{2} \left( \sin \frac{W}{2} \cos \frac{X}{2} - \cos \frac{W}{2} \sin \frac{X}{2} \right) &  \\
  + 2 \sin \frac{W}{2} \cos ^{2} \frac{X}{2} \sin \frac{Y}{2}\left( \sin \frac{W'}{2} \cos \frac{Y}{2} - \cos \frac{W'}{2} \sin \frac{Y}{2} \right) &  
\end{array}$

\vspace{2 mm}

\noindent Note that the last two quatities in parentheses are nothing but $\displaystyle{\sin (W-X)/2}$ and  $\displaystyle{\sin (W'-Y)/2}$. Combined with the first two rows gives us

\vspace{2 mm}

$  2 \cos \frac{X}{2} \sin \frac{Y}{2} \sin \left( \frac{W' - Y}{2} \right) \left( - \sin \frac{X}{2} \cos \frac{W}{2} + \sin \frac{W}{2} \cos \frac{X}{2} \right)$

\vspace{2 mm}

$+  2 \sin \frac{X}{2} \cos \frac{Y}{2} \sin   ( \frac{W-X}{2} ) \left( - \sin \frac{Y}{2} \cos \frac{W'}{2} + \cos \frac{Y}{2} \sin \frac{W'}{2} \right)$

\vspace{2 mm}

$ = 2  \cos \frac{X}{2} \sin \frac{Y}{2} \sin \left( \frac{W' - Y}{2} \right) \sin \left( \frac{W-X}{2} \right) + 2 \sin \frac{X}{2} \cos \frac{Y}{2} \sin   \left( \frac{W-X}{2} \right) \sin \left( \frac{W' - Y}{2} \right)$

\vspace{2 mm} 

$=  2   \sin \left( \frac{W' - Y}{2} \right) \sin \left( \frac{W-X}{2} \right) \sin \left( \frac{X+Y}{2} \right) \, \geq 0$

\vspace{2 mm}

\noindent This last expression is nonnegative because $W'-Y= \alpha _{3} + \beta _{2}$, $W-X= \alpha _{4} + \beta _{3}$, and  these angles will be between $0$ and $\pi$ when divided by 2 (in addition to  $X+Y$), which conlcudes the proof of Theorem 2.

\vspace{2 mm}

\underline{Remark}: Now that Theorem 2 has been established, let us go back and collect all remaining terms on the left hand side of (4). We obviously have the final answer which is $2   \sin \frac{W' - Y}{2} \sin  \frac{W-X}{2}  \sin  \frac{X+Y}{2} $. However, we have some additional terms remaining from corollary 1. Namely, we have the following three terms from the inequalities used in Corollary 1 (substituted in the left column of Theorem 2) which remain on the left hand side of (4):

\vspace{2 mm}

$ \, \, \, \, \,  2 \sin X  \sin \frac{W'}{2}  \sin \frac{Y}{2} \sin  \frac{\alpha _{3}}{2} \cos \frac{\beta _{2}}{2}  $

$+ 2 \sin Y  \sin \frac{W}{2}  \sin \frac{X}{2} \sin \frac{\beta _{3}}{2} \cos \frac{\alpha _{4}}{2}    $  \hfill (5)

$+ 2  \sin W  \cos \frac{X}{2}  \cos \frac{Y}{2} \cos  \frac{\gamma _{1}}{2} \sin \frac{\gamma _{3}}{2} $

\vspace{2 mm}

In order to finish the proof of the quadilaterial inequality (1), we need to show that the remaining terms give a positive answer when added to the second line of Lemma 2, which is $ - 2 \sin  \frac{\beta _{4} -\alpha _{1}}{2} \sin  \frac{\alpha _{2} -\beta _{1}}{2} \sin \frac{\gamma _{3} + \gamma _{1}}{2}$. Well, if the quadrilateral at hand satisfies $\sin \frac{\beta _{4} -\alpha _{1}}{2} \sin  \frac{\alpha _{2} -\beta _{1}}{2} \leq 0$, then our result is immediate. Accordingly, we still need to prove our inequality (1) under the  assumption that $\sin \frac{\beta _{4} -\alpha _{1}}{2} \sin  \frac{\alpha _{2} -\beta _{1}}{2} \geq 0$. Let us start by adding the resulting positive term $2 \sin  \frac{W' - Y}{2}  \sin  \frac{W-X}{2}  \sin  \frac{X+Y}{2}$ from Theorem 2 to the second line of Lemma 2:

\vspace{2 mm} 

$\displaystyle{2 \sin  \frac{W' - Y}{2}  \sin  \frac{W-X}{2}  \sin  \frac{X+Y}{2}  -  2 \sin  \frac{\beta _{4} -\alpha _{1}}{2} \sin  \frac{\alpha _{2} -\beta _{1}}{2} \sin \frac{\gamma _{3} + \gamma _{1}}{2}}$

$\displaystyle{ \geq 2 \sin  \frac{\beta _{4} -\alpha _{1}}{2} \sin  \frac{\alpha _{2} -\beta _{1}}{2} \left(  \sin  \frac{X+Y}{2} - \sin \frac{\gamma _{3} + \gamma _{1}}{2}\right)}$

$\displaystyle{ = - 4 \sin  \frac{\beta _{4} -\alpha _{1}}{2} \sin  \frac{\alpha _{2} -\beta _{1}}{2}  \cos \frac{\gamma _{1}}{2} \sin \frac{\gamma _{3}}{2}}$

\vspace{2 mm}

\noindent where the first inequality follows from Corollary 1 and the last equality comes from the fact that $X+Y=\gamma _{1}- \gamma _{3}$. We can see that, under our assumption  $\sin \frac{\beta _{4} -\alpha _{1}}{2} \sin  \frac{\alpha _{2} -\beta _{1}}{2} \geq 0$, our outcome is still negative. Accordingly, we take this answer and add to it the third line of (5) to get 

\vspace{2 mm}

$\displaystyle{ 4 \sin \frac{W}{2} \cos \frac{W}{2}  \cos \frac{X}{2} \cos \frac{Y}{2} \cos \frac{\gamma _{1}}{2} \sin \frac{\gamma _{3}}{2}}$

\vspace{2 mm} 

$\displaystyle{  - 4 \sin  \frac{\beta _{4} -\alpha _{1}}{2} \sin  \frac{\alpha _{2} -\beta _{1}}{2}  \cos \frac{\gamma _{1}}{2} \sin \frac{\gamma _{3}}{2}}$

\vspace{2 mm}

This is easily seen to be positive: We factor out $4 \cos \frac{\gamma _{1}}{2} \sin \frac{\gamma _{3}}{2}$ which is positive. We then note that $\sin \frac{W}{2} \cos \frac{X}{2} \geq \sin \frac{W-X}{2}$ and $\sin \frac{W'}{2} \cos \frac{Y}{2} \geq \sin \frac{W'-Y}{2}$, which, by Corollary (1), are not less than $\sin  \frac{\alpha _{1} - \beta _{4} }{2}$ and $ \sin  \frac{\beta _{1} - \alpha _{2} }{2}$ respectively. Since $\sin  \frac{\beta _{4} -\alpha _{1}}{2} \sin  \frac{\alpha _{2} -\beta _{1}}{2} \geq 0$, we can see that we end up with a non-negative answer,  which finally finishes the proof of inequality (1).


\begin{thebibliography}{20}
%
\bibitem{A1} M.F. Atiyah,   \emph{The geometry of classical particles}. Surveys in differential geometry $\bf{vol VII}$, 1-15, International Press, Somerville, MA  (2000)
%
\bibitem{A2}  M.F. Atiyah: \emph{Configuration of points}.  Phil. Trans. R. Lond.  $\bf{A 359}$, 1375-1387 (2001)
%
\bibitem{AS} Michael Atiyah and Paul Sutcliffe, The geometry of point particles, R. Soc. Lond. Proc. Ser. A Math.
Phys. Eng. Sci. 458 (2002), no. 2021, 1089–1115. MR 1902577 (2003c:55019)
%
\bibitem{BR}  M.V. Berry, J.M. Robbins: \emph{Indistinguishability for quantum particles: spin, statistics and geometric phase}. Proc. Roy. Soc. London Ser. . $\bf{A 453}$,  1771-1790 (1997)
%
\bibitem{D}  D. Z. Dokovic: \emph{Verification of Atiyah's conjecture for some nonplanar configurations with dihedral symmetry}. Publ. Inst. Math. (Beograd) (N.S.) $\bf{72}$, 23-28 (2002)
%
\bibitem{EN} M. Eastwood, P. Norbury: \emph{A proof of Atiyah's conjecture on configurations of four points in Euclidean three-space}. Geom. Topol. $\bf{5}$, 885-893  (2001)
%
\bibitem{M} Malkoun, J. Configurations of points and the symplectic Berry-Robbins problem,
SIGMA Symmetry Integrability Geom. Methods Appl. 10 (2014), Paper 112,
6pp.
%
\bibitem{MP} M. Mazur, B. V. Penterko: \emph{On the conjectures of Atiyah and Sutcliffe}. Geom. Dedicata. DOI 10.1007/s10711-011-9636-6




\end{thebibliography}
\end{document}